\documentclass[11pt,a4paper]{article}

\usepackage{latexsym,amssymb, amsthm}
\usepackage[centertags]{amsmath}
\usepackage{epsfig,pifont}
\usepackage[inactive]{srcltx}  

\usepackage{graphics,helvet,times,mathptm,epic,eepic}

\usepackage{color}

\usepackage{newlfont}

\usepackage{amsfonts,bbm}

\usepackage{url}

\newenvironment{keywords}{ \noindent {\small\bf Key Words}:}{ }

\def\bd{\begin{description}}
\def\ed{\end{description}}

\def\beq{\begin{equation}}
\def\eeq{\end{equation}}

\def\G1{\hbox{$\displaystyle{\mbox{\ding{172}}}$}}

\def\bea{\begin{eqnarray}}
\def\eea{\end{eqnarray}}
\def\beas{\begin{eqnarray*}}
\def\eeas{\end{eqnarray*}}

\newtheorem{theorem}{Theorem}[section]

\theoremstyle{remark}

\begin{document}

\title{\textbf{On accuracy of mathematical languages used\\ to deal with
the Riemann zeta function and\\ the Dirichlet eta function}}

\newcommand{\nms}{\normalsize}
\author{  \bf Yaroslav D. Sergeyev\footnote{Yaroslav D. Sergeyev, Ph.D., D.Sc., is Distinguished
Professor at the University of Calabria, Rende, Italy.
 He is also Full Professor   at the N.I.~Lobatchevsky State University,
  Nizhni Novgorod, Russia  and Affiliated Researcher at the Institute of High Performance
Computing and Networking of the National Research Council of
Italy.
 {\tt yaro@si.deis.unical.it}, \,\, {\tt
http://wwwinfo.deis.unical.it/$\sim$yaro} }  }
\date{}

\maketitle

\vspace{-1cm}

\begin{abstract}

The Riemann Hypothesis  has been of central interest to
mathematicians for a long time  and many unsuccessful attempts
have been made to either prove or disprove it. Since the Riemann
zeta function is defined as a sum of the infinite number of items,
in this paper, we look at the Riemann Hypothesis using a new
applied approach to infinity allowing one to easily execute
numerical computations with various infinite and infinitesimal
numbers  in accordance with the principle `The part is less than
the whole' observed in the physical world around us.    The new
approach allows one to work with functions and derivatives   that
can assume not only finite but also infinite  and infinitesimal
values and this possibility is used to study properties of  the
Riemann zeta function and the Dirichlet eta function. A new
computational approach allowing one to evaluate these functions at
certain points is proposed. Numerical examples are given. It is
emphasized that different mathematical languages can be used to
describe mathematical objects with different accuracies. The
traditional and the new approaches are compared with respect to
their application to the Riemann zeta function and the Dirichlet
eta function. The accuracy of the obtained results is discussed in
detail.
 \end{abstract}

\begin{keywords}
Infinite and infinitesimal numbers and numerals; accuracy of
mathematical languages; Riemann zeta function; Dirichlet eta
function; divergent series.
 \end{keywords}




\section{Introduction}
\label{s1_riem}
 The Riemann zeta function
 \beq
   \zeta(s)=\sum_{u=1}^{\infty} \frac{1}{u^s}
  \label{Riemann_zeta}
 \eeq
defined for complex numbers $s$ is one of the most important
mathematical objects discovered so far (see
\cite{Conrey,Edwards}).  It has been first introduced and studied
by Euler (see \cite{Euler_1,Euler_2,Euler_3,Euler_4}) who proved
the famous identity
 \beq
   \sum_{u=1}^{\infty} \frac{1}{u^s} =
    \prod_{p\,\, \mbox{\tiny{primes}}}  \frac{1}{1-p^{-s}}
  \label{Riemann_primes}
 \eeq
being one of the strongest sources of the interest to the Riemann
zeta function. Many interesting results have been established for
this function (see \cite{Conrey,Edwards} for a complete
discussion). In the context of this paper, the following results
will be of our primary interest.

It is  known that $\zeta(s)$ diverges on the open half-plane of
$s$ such that the real part $\Re(s) < 1$. Then, for $s=0$ the
following value is attributed to $\zeta(0)$:
 \beq
 \zeta(0)=\sum_{u=1}^{\infty} \frac{1}{u^0}=  1+1+1+\ldots  =
 -\frac{1}{2}.
  \label{1+1}
 \eeq
It is also known that the Riemann zeta function has the following
relation
 \beq
 \eta(s) = (1- 2^{1-s})\zeta(s)
   \label{Riemann_relation}
 \eeq
to the Dirichlet eta function
 \beq
   \eta(s)=\sum_{u=1}^{\infty} \frac{{(-1)}^{u-1}}{u^s}.
  \label{Riemann_eta}
 \eeq
 It has been shown for   the Riemann
zeta function that it has trivial  zeros at the points
$-2,-4,\dots$  It is also known that any non-trivial zero lies in
the complex set $\{s : 0 < \Re(s) < 1 \}$,    called the critical
strip. The complex set $\{s :   \Re(s) =1/2 \}$ is called the
critical line. The Riemann Hypothesis asserts that any non-trivial
zero $s$ has the real part $\Re(s) = 1/2$, i.e., lies on the
critical line.

This problem has been of central interest to mathematicians for a
long time  and many unsuccessful attempts have been made to either
prove or disprove it. Since both the Riemann zeta function and the
Dirichlet eta function are defined as   sums of the infinite
number of items, in this paper, we look at them and  at the
Riemann Hypothesis using a new approach to infinity introduced in
\cite{Sergeyev,www,informatica,first,Lagrange}. This approach, in
particular, incorporates the following two ideas.

i) The point of view on infinite and infinitesimal quantities
applied in this paper uses strongly two methodological ideas
borrowed from the modern Physics: relativity and interrelations
holding between the object of an observation and the tool used for
this observation. The latter is directly related to interrelations
between numeral systems\footnote{ We remind that \textit{numeral}
is a symbol or group of symbols that represents a \textit{number}.
The difference between numerals and numbers is the same as the
difference between words and the things they refer to. A
\textit{number} is a concept that a \textit{numeral} expresses.
The same number can be represented by different numerals. For
example, the symbols `3', `three', and `III' are different
numerals, but they all represent the same number.} used to
describe numbers and numbers (and other mathematical objects)
themselves. Numerals that we use to write down numbers, functions,
etc. are among our tools of investigation and, as a result, they
strongly influence our capabilities to study mathematical objects.

ii) Both standard and non-standard Analysis study mainly functions
assuming finite values. In \cite{informatica,Dif_Calculus},
functions and their derivatives that can assume finite, infinite,
and infinitesimal values and can be defined over finite, infinite,
and infinitesimal domains have been studied. Infinite and
infinitesimal numbers are not auxiliary entities in the new
approach, they are full members in it and can be used in the same
way as finite constants.

In the next section, we present very briefly the new methodology
for treating infinite and infinitesimal quantities indicating
mainly the facts that are used directly in this paper. A
comprehensive introduction and numerous examples of its usage can
be found in \cite{informatica,Lagrange} while some applications
can be found in \cite{De_Leone,D'Alotto,Margenstern,informatica,
Korea,Dif_Calculus,Num_dif,Sergeyev_Garro}.

\section{Accuracy of mathematical languages and a new numeral system for dealing with   infinity}
\label{s2_riem}

In order to introduce the new methodology, let us consider a study
published in \textit{Science} by Peter Gordon (see \cite{Gordon})
where he describes a primitive tribe living in Amazonia --
Pirah\~{a} -- that uses a very simple numeral system  for
counting: one, two, many. For Pirah\~{a}, all quantities larger
than two are just `many' and such operations as 2+2 and 2+1 give
the same result, i.e., `many'. Using their weak numeral system
Pirah\~{a} are not able to see, for instance, numbers 3, 4,  and
5, to execute arithmetical operations with them, and, in general,
to say anything about these numbers because in their language
there are neither words nor concepts for that.

It is important to emphasize that the records $ 2+1=\mbox{`many'}$
and $2+2= \mbox{`many'}$ are correct in their language and if one
is satisfied with the accuracy of  the answer `many', it can be
used (and \textit{is used} by Pirah\~{a}) in practice. Note that
also for us, people knowing that $2+1=3$ and $2+2=4$, the result
of Pirah\~{a} is not wrong, it is just \textit{inaccurate}. Thus,
if one needs a more precise result than `many', it is necessary to
introduce a more powerful mathematical language (a numeral system
in this case) allowing one to express the required answer in a
more accurate way. By using modern powerful numeral systems where
additional numerals  for expressing the numbers `three' and `four'
have been introduced, we can notice that within `many' there are
several objects and the numbers 3 and 4 are among these unknown to
Pirah\~{a}) objects.

Thus, the choice of the mathematical language depends on the
practical problem that is to be   solved and on the accuracy
required for such a solution. In dependence of this accuracy, a
numeral system  that would be able to express the numbers
composing the answer should be chosen.

Such a situation  is typical for   natural sciences where it is
well known that instruments bound and influence results of
observations. When physicists see a black dot in their microscope
they cannot say: The object of the observation \textit{is} the
black dot. They are obliged to say: the lens used in the
microscope allows us to see the black dot and it is not possible
to say anything more about the nature of the object of the
observation until we  change the instrument - the lens or the
microscope itself - by a more precise one. Then, probably, when a
stronger lens will be put in the microscope, physicists will be
able to see that the object that seemed to be one dot  with the
new lens can be viewed as, e.g., two dots.

Note that both results (one dot and two dots) correctly represent
the reality with the accuracy of the chosen instrument of the
observation (lens). Physicists decide the level of the precision
they need and obtain a result depending of the chosen level of the
accuracy. In the moment when they put a lens in the microscope,
they have decided the minimal (and the maximal) size of objects
that they will be able to observe. If they need a more precise or
a more rough answer, they change the lens of their microscope.

Analogously, when mathematicians have decided which mathematical
languages (in particular, which numeral systems) they will use in
their work, they have decided which mathematical objects they will
be able to observe and to describe.

In natural sciences   there always exists the triad -- the
researcher, the object of investigation, and tools used to observe
the object -- and  the instrument used to observe the object
bounds and influences results of observations. The same happens in
Mathematics  studying numbers and objects that can be constructed
by using numbers. Numeral systems used to express numbers are
instruments of observations used by mathematicians. The usage of
powerful numeral systems gives the possibility of obtaining more
precise results in Mathematics in the same way as usage of a good
microscope gives the possibility to obtain more precise results in
Physics. Let us now return to Pirah\~{a} and make the first
relevant observation with respect to the subject of this paper.

The weakness of the  numeral system of Pirah\~{a}   gives them
inaccurate answers to arithmetical operations with finite numbers
where modern numeral systems are able to provide exact solutions.
However, Pirah\~{a}  obtain also such results as
 \beq
\mbox{`many'}+ 1= \mbox{`many'},   \hspace{1cm}    \mbox{`many'} +
2 = \mbox{`many'},
    \label{Riemann29}
 \eeq
which are very familiar to us  in the context of views on infinity
used in the traditional calculus where the famous numeral $\infty$
is used:
 \beq
\infty + 1= \infty,    \hspace{1cm}    \infty + 2 = \infty.
    \label{Riemann30}
 \eeq
This observation leads us to the following idea: \textit{Probably
our difficulties in working with infinity in general (and in
particular, with the the Riemann zeta function) is not connected
to the nature of infinity itself but is a result of inadequate
numeral systems that we use to work with infinity, more precisely,
to express infinite numbers.}

Recently,   a new   numeral system has been developed in order to
express finite, infinite, and infinitesimal numbers in a unique
framework (a rather comprehensive description of the new
methodology  can be found  in \cite{informatica,Lagrange}). The
main idea consists of measuring infinite and infinitesimal
quantities by different (infinite, finite, and infinitesimal)
units of measure. This is done by using a new numeral.

An infinite unit of measure   has been introduced  for this
purpose in \cite{informatica,Lagrange} as the number of elements
of the set $\mathbb{N}$ of natural numbers. It is expressed by a
new numeral \G1 called \textit{grossone}. It is necessary to
emphasize immediately that the infinite number \G1 is neither
Cantor's $\aleph_0$ nor $\omega$ and the new approach is not
related to the non-standard analysis.

One of the important differences with respect to traditional views
consists of the fact that the new approach continuously emphasizes
the distinction between numbers and numerals and studies in depth
numeral systems. In general, neither standard analysis nor
non-standard one pay a lot of attention to the necessity, in order
to be able to execute an arithmetical operation, to have on hand
numerals required to express both the operands and the result of
the operation. The accuracy of the obtained result (i.e., how well
the numeral chosen to represent the result of the operation
expresses the resulting quantity) is not studied either. For the
approach introduced in \cite{informatica,Lagrange} both issues are
crucial. Since for any fixed numeral system $\mathcal{S}$ there
exist operations that cannot be executed in $\mathcal{S}$ because
$\mathcal{S}$ has no suitable numerals, these marginal situations
are investigated in detail. We can remind Pirah\~{a} in this
occasion (see (\ref{Riemann29}), (\ref{Riemann30})) but even the
modern powerful positional numeral systems which are not able to
work with numerals having, e.g., $10^{100}$ digits. It is shown in
\cite{informatica,first,Lagrange,Sergeyev_Garro} that very often
the margins of the expressibility of numeral systems are connected
directly to interesting research problems.

Another important peculiarity of the new approach consists of the
fact that \G1 has both cardinal and ordinal properties as usual
finite natural numbers have. In fact, infinite positive integers
that can be viewed through numerals including grossone  can be
interpreted in the terms of the number of elements of certain
infinite sets. For instance, the set of even numbers has
$\frac{\G1}{2}$ elements and the set of integers has
$2\mbox{\G1\small{+}}1$ elements. Thus, the new numeral system
allows one to distinguish within countable sets many different
sets having the different infinite number  of elements (remember
again the analogy with the microscope). Analogously, within
uncountable sets it is possible to distinguish sets having, for
instance, $2^{\mbox{\tiny{\G1}}}$ elements,
$10^{\mbox{\tiny{\G1}}}$ elements,   and even
$\G1^{\mbox{\tiny{\G1}}}-1$, $\G1^{\mbox{\tiny{\G1}}}$, and
$\G1^{\mbox{\tiny{\G1}}}+1$ elements and to show  (see
\cite{first,Lagrange}) that
 \[
\frac{\G1}{2} < \G1 < 2\mbox{\G1\small{+}}1 <
2^{\mbox{\tiny{\G1}}} < 10^{\mbox{\tiny{\G1}}} <
\G1^{\mbox{\tiny{\G1}}}-1 < \G1^{\mbox{\tiny{\G1}}}<
\G1^{\mbox{\tiny{\G1}}}+1.
\]

Formally, grossone is introduced as a new number by describing its
properties postulated by the \textit{Infinite Unit Axiom} (see
\cite{informatica,Lagrange}). This axiom is added to axioms for
real numbers similarly to addition of the axiom determining zero
to axioms of natural numbers when integers are introduced.
Inasmuch as it has been postulated that grossone is a number,  all
other axioms for numbers hold for it, too. Particularly,
associative and commutative properties of multiplication and
addition, distributive property of multiplication over addition,
existence of   inverse  elements with respect to addition and
multiplication hold for grossone as for finite numbers. This means
that  the following relations hold for grossone, as for any other
number
 \beq
 0 \cdot \G1 =
\G1 \cdot 0 = 0, \hspace{3mm} \G1-\G1= 0,\hspace{3mm}
\frac{\G1}{\G1}=1, \hspace{3mm} \G1^0=1, \hspace{3mm}
1^{\mbox{\tiny{\G1}}}=1, \hspace{3mm} 0^{\mbox{\tiny{\G1}}}=0.
 \label{3.2.1}
       \eeq

   Various numeral
systems including \G1 can be used for expressing infinite and
infinitesimal numbers. In particular, records similar to
traditional positional numeral systems can be used (see
\cite{informatica,Korea}).
In order to construct a number
$C$ in this system, we subdivide $C$ into groups corresponding to
powers of grossone:
 \beq
  C = c_{p_{m}}
\G1^{p_{m}} +  \ldots + c_{p_{1}} \G1^{p_{1}} +c_{p_{0}}
\G1^{p_{0}} + c_{p_{-1}} \G1^{p_{-1}} + \ldots   + c_{p_{-k}}
 \G1^{p_{-k}}.
\label{3.12}
       \eeq
 Then, the record
 \beq
  C = c_{p_{m}}
\G1^{p_{m}}    \ldots   c_{p_{1}} \G1^{p_{1}} c_{p_{0}}
\G1^{p_{0}} c_{p_{-1}} \G1^{p_{-1}} \ldots c_{p_{-k}}
 \G1^{p_{-k}}
 \label{3.13}
       \eeq
represents  the number $C$, where finite numbers $c_i\neq0$ called
\textit{grossdigits}   can be   positive or negative. They show
how many corresponding units should be added or subtracted in
order to form the number $C$. Grossdigits can be expressed by
several symbols.

Numbers $p_i$ in (\ref{3.13}) called \textit{grosspowers}  can be
finite, infinite, and infinitesimal, they   are sorted in the
decreasing order with $ p_0=0$
\[
p_{m} >  p_{m-1}  > \ldots    > p_{1} > p_0 > p_{-1}  > \ldots
p_{-(k-1)}  >   p_{-k}.
 \]
 In the record (\ref{3.13}), we write
$\G1^{p_{i}}$ explicitly because in the new numeral positional
system  the number   $i$ in general is not equal to the grosspower
$p_{i}$.

 Finite numbers  in this new numeral system are represented
by numerals having only one grosspower $ p_0=0$. In fact, if we
have a number $C$ such that $m=k=$~0 in representation
(\ref{3.13}), then due to (\ref{3.2.1}),   we have $C=c_0
\G1^0=c_0$. Thus, the number $C$ in this case does not contain
grossone and is equal to the grossdigit $c_0$ being a conventional
finite number     expressed in a traditional finite numeral
system.

 The simplest infinitesimal numbers    are represented by numerals $C$
having only finite or infinite negative  grosspowers, e.g.,
$21.3\G1^{-243.6}5\G1^{-154\mbox{\tiny{\G1}}}$. Then, the
infinitesimal number $\frac{1}{\G1}=\G1^{-1}$   is the inverse
element with respect to multiplication for~\G1:
 \beq
\G1^{-1}\cdot\G1=\G1\cdot\G1^{-1}=1.
 \label{3.15.1}
       \eeq
Note that all infinitesimals are not equal to zero. Particularly,
$\frac{1}{\G1}>0$ because it is a result of division of two
positive numbers.

The simplest infinite numbers     are expressed by numerals having
  positive finite or infinite  gross\-powers without infinitesimals. They have
infinite parts and can also have  a finite part and infinitesimal
ones.  For instance, the number
\[
 31.5\G1^{14.8\mbox{\tiny{\G1}}}(-0.645)\G1^{5}
 7.89\G1^{0}37\G1^{-4.29}
 72.8\G1^{-360.21}
\]
has two infinite parts $31.5\G1^{14.8\mbox{\tiny{\G1}}}$ and
$-0.645\G1^{5}$ one finite part $7.89\G1^{0}$ and two
infinitesimal parts $37\G1^{-4.29}$ and $72.8\G1^{-360.21}$. All
of the numbers introduced above can be used as grosspowers, as
well, giving so a possibility to have various combinations of
quantities and to construct  terms having a more complex
structure. In this paper, we are interested only on numbers having
integer grosspowers.

\section{Sums having a fixed infinite number  of items and\\ integrals
 over fixed infinite domains}
\label{s3_riem}

Introduction of the new numeral, \G1, allows us to distinguish
within $\infty$ many different infinite numbers in the same way as
the introduction of modern numeral systems used to express finite
numbers allows us to distinguish numbers `three' and `four' within
`many'.

Thus, the Aristotle principle `The part is less than the whole'
leading to the fact that $x+1>x$ can be applied to all numbers $x$
(finite, infinite, and infinitesimal) and  not only to finite
numbers as it is usually done in the traditional mathematics (this
point is discussed in detail in \cite{informatica,Lagrange}). In
general, we   have that $x + y> x,\,\, y>0,$ where both $y$ and
$x$ can be
  finite, infinite, or infinitesimal. For instance,
it follows
$$\G1+1 > \G1, \hspace{1cm} 1+\G1^{-1} > 1,
\hspace{1cm}\G1^2+1 > \G1^2, \hspace{1cm} \G1^2+\G1^{-1} >
\G1^2. $$

 This means
that such records as $S=a_1+a_2+\ldots$ or
$\sum_{i=1}^{\infty}a_i$ become unprecise in the new mathematical
language using grossone. By continuation the analogy with
Pirah\~{a} (see     (\ref{Riemann29}), (\ref{Riemann30})), the
record $\sum_{i=1}^{\infty}a_i$ becomes a kind of
$\sum_{i=1}^{many}a_i$.

It is worthwhile noticing also that the symbol $\infty$ cannot be
used in the same expression with numerals using \G1 (analogously,
the record `many'+4 has no sense because it uses symbols from two
different languages having different accuracies).

 As a consequence, if one wants to work
with an infinite number, $n$, of items in a sum $\sum_{i=0}^{n}
a_i$ then it is necessary to fix explicitly this number using
numerals available  for expressing infinite numbers in a chosen
numeral system (e.g., the numeral system introduced in the
previous section can be taken  for this purpose). The same
situation we have with sums having a finite number of items: it is
not sufficient to say that the number, $n$, of items in the sum is
finite, it is necessary to fix explicitly the value of~$n$ using
for this purpose numerals available in a chosen traditional
numeral system to express finite numbers.

The new numeral system using grossone allows us to work with
expressions involving   infinite numbers   and to obtain, where it
is possible, as results infinite, finite, and infinitesimal
numbers. Obviously, this is of a particular interest in connection
with the Riemann zeta function.  For instance,  it becomes
possible to reconsider  the result (\ref{1+1}) that is very
difficult to be fully realized by anyone who is not a
mathematician. In fact, when one has an infinite sum of positive
finite integers he or she would expect to have an infinite
positive integer as a result. In contrast, (\ref{1+1}) proposes a
negative finite fractional number as the answer.

In order to become acquainted with the way   the new methodology
is applied to the theory of divergent series, let us consider
several example. We start by considering   two infinite series
$S_1=7+7+7+7+\ldots$ and $S_2=3+3+3+3+\ldots$   The traditional
analysis gives us a very poor answer that both of them diverge to
infinity. Such operations as, e.g., $\frac{S_2}{S_1}$ and $S_2 -
S_1$
 are not defined.

The new mathematical language using grossone allows us to indicate
explicitly the number of their items. Suppose that the   series
$S_1$ has $k$ items and $S_2$ has $n$ items:
$$S_1(k)=\underbrace{7+7+7+\ldots+7}_{k \mbox{ items}},
\hspace{1cm} S_2(n)=\underbrace{3+3+3+\ldots+3}_{n \mbox{
items}}.$$ Then $S_1(k)=7k$ and $S_2(n)=3n$ and by giving
different numerical values (finite or infinite) to $k$ and $n$ we
obtain different numerical values for the sums.  For chosen $k$
and $n$ it becomes possible to calculate $S_2(n) - S_1(k)$
(analogously, the expression $ \frac{S_1(k)}{S_2(n)} $ can be
calculated). If, for instance, $k=5\G1$ and $n=\G1$   we obtain
$S_1(5\G1)=35\G1$, $S_2(\G1)=3\G1$ and it follows
\[
S_2(\G1) - S_1(5\G1)= 3\G1 - 35\G1 = - 32\G1< 0.
\]
If  $k=3\G1$ and $n=7\G1+2$   we obtain $S_1(3\G1)=21\G1$,\,\,
$S_2(7\G1+2)=21\G1+6$ and  it follows
\[
S_2(7\G1+2) -  S_1(3\G1)= 21\G1+6 - 21\G1 = 6.
\]
 It becomes
also possible to calculate sums having an infinite number of
infinite or infinitesimal items. Let us consider, for instance,
sums
\[
S_3(l)=\underbrace{2\G1+2\G1+\ldots+2\G1}_{l \mbox{ items}},
\hspace{1cm}
S_4(m)=\underbrace{4\G1^{-1}+4\G1^{-1}+\ldots+4\G1^{-1}}_{m \mbox{
items}}.
\]
For $l=m=0.5\G1$ it follows $S_3(0.5\G1)= \G1^2$ and
$S_4(0.5\G1)=2$ (remind that $\G1\cdot\G1^{-1}=\G1^{0}=1$ (see
(\ref{3.15.1})). It can be seen from this example that it is
possible to obtain finite   numbers as the result of summing up
infinitesimals.

The infinite and infinitesimal numbers allow us to calculate also
arithmetic and geometric series with an infinite number of items.
For the arithmetic  progression,  $a_{n} = a_{1} + (n - 1)d$, for
both finite and infinite  $n$ we have
\[
\sum_{i=1}^{n} a_{i} = \frac{n}{2}(a_{1} + a_{n}).
\]
Then, for instance, the sum of all natural numbers from 1 to \G1
can be calculated as follows
 \beq
1+2+3+ \ldots + (\G1-1) + \G1 = \sum_{i=1}^{\G1} i =
\frac{\G1}{2}(1 + \G1)= 0.5\G1^{2}0.5\G1.
 \label{Riemann1}
 \eeq
 If we calculate  now the following sum of
infinitesimals where each item is \G1 times less than the
corresponding item of (\ref{Riemann1})
\[
\G1^{-1}+2\G1^{-1}+ \ldots + (\G1-1)\cdot\G1^{-1} +
\G1\cdot\G1^{-1} = \sum_{i=1}^{\G1} i\G1^{-1} =
\frac{\G1}{2}(\G1^{-1} + 1)= 0.5\G1^{1}0.5.
\]
then, obviously, the obtained number, $0.5\G1^{1}0.5$, is \G1
times less than the sum in (\ref{Riemann1}).  This example shows,
particularly, that infinite numbers can also be obtained as the
result of summing up infinitesimals.

Let us consider now the geometric series $\sum_{i=0}^{\infty}
x^i$. Traditional analysis proves that it converges to
$\frac{1}{1-x}$ for $x$ such that $-1 < x < 1$. We are able to
give a more precise answer for \textit{all} values of $x$. To do
this, we should fix the number of items in the sum. If we suppose
that it contains $n$ items, where $n$ is finite or infinite, then
 \beq
 Q_n = \sum_{i=0}^{n}
x^i = 1 + x + x^2 + \ldots + x^n.
 \label{3.7.2.f}
 \eeq
By multiplying the left hand and the right hand parts of this
equality by $x$ and by subtracting the result from (\ref{3.7.2.f})
we obtain
\[
Q_n - xQ_n = 1-x^{n+1}
\]
and, as a consequence, for all $x\neq 1$ the formula
 \beq
  Q_n =
(1-x^{n+1})(1-x)^{-1}
 \label{3.7.2.f.1}
 \eeq
holds for finite and infinite $n$. Thus, the possibility to
express infinite and infinitesimal numbers allows us   to take
into account    infinite $n$   and the value $x^{n+1}$ being
infinitesimal for a finite $x$. Moreover, we can calculate $Q_n$
for infinite and finite values of $n$ and $x=1$, because in this
case we have just
\[
Q_n = \underbrace{1+1+1+\ldots+1}_{n+1 \mbox{ items}} = n+1.
\]

Let us illustrate the obtain results by some examples. In the
first of them we consider the divergent series
\[
1 + 3 + 9 + \ldots = \sum_{i=0}^{\infty} 3^i.
\]
In the new language, to define a sum it is necessary to fix the
number of items, $n$, in it and $n$ can be finite or infinite. For
example, for the infinite $n=\G1^2$ we obtain
\[
\sum_{i=0}^{\G1^2} 3^i = 1 + 3 + 9 + \ldots +
3^{\mbox{\tiny{\G1}}^2}= \frac{1-3^{\mbox{\tiny{\G1}}^2+1}}{1-3} =
0.5(3^{\mbox{\tiny{\G1}}^2+1}-1).
\]
Analogously, for  $n=\G1^2+1$ we obtain
\[
1 + 3 + 9 + \ldots + 3^{\mbox{\tiny{\G1}}^2} +
3^{\mbox{\tiny{\G1}}^{2}+1}= 0.5(3^{\mbox{\tiny{\G1}}^2+2}-1).
\]
If we now find the difference between the two sums
\[
0.5(3^{\mbox{\tiny{\G1}}^2+2}-1) -0.5(3^{\mbox{\tiny{\G1}}^2+1}-1)
= 3^{\mbox{\tiny{\G1}}^2+1} (0.5\cdot3 - 0.5) =
3^{\mbox{\tiny{\G1}}^{2}+1}
\]
 we  obtain the
newly added item $3^{\mbox{\tiny{\G1}}^2+1}$.

Let us now consider the series $\sum_{i=1}^{\infty}\frac{1}{2^i}$.
It is well known that it converges to one. However, we are able to
give a more precise answer. In fact, the formula
 \[
 \sum_{i=1}^{n}\frac{1}{2^i}= \frac{1}{2}( 1+ \frac{1}{2} + \frac{1}{2^2} + \ldots + \frac{1}{2^{n-1}}) = \frac{1}{2} \cdot \frac{1-\frac{1}{2}^{n}}{1-\frac{1}{2}}  = 1-\frac{1}{2^n}
 \]
can be used directly for   infinite $n$, too. For example, if
$n=\G1$ then
$$\sum_{i=1}^{\mbox{\small{\G1}}}\frac{1}{2^i}=1-\frac{1}{2^{\mbox{\tiny{\G1}}}},$$
where $\frac{1}{2^{\mbox{\tiny{\G1}}}}$ is infinitesimal. Thus,
the traditional answer $\sum_{i=1}^{\infty}\frac{1}{2^i}=1$ was
just a finite approximation to our more precise result using
infinitesimals.

In order to show the potential of the new approach for the work
with divergent series with alternate signs, we start from the
famous series
\[
S_5 = 1-1+1-1+1-1+\ldots
\]
 In literature, there exist  many
approaches giving different answers regarding the value of this
series (see \cite{Knopp}). All of them use various notions of
average to calculate the series. However, the notions of the sum
and an average are different. In our approach, we do not appeal to
an average and calculate the required sum directly. To do this we
should indicate explicitly the number of items, $k$, in the sum.
Then
\beq
S_5(k)=\underbrace{1-1+1-1+1-1+1-\ldots}_{k \mbox{
items}} = \left \{
\begin{array}{ll} 0, &
  \mbox{if  } k=2n,\\
1, &    \mbox{if  } k=2n+1,\\
 \end{array} \right.
 \label{Riemann27}
\eeq
 where $k$ can be finite or infinite. For
example, $S_5(\G1)=0$ because   \G1 is even (see
\cite{informatica}). Analogously, $S_5(\G1-1)=S_5(2\G1+1)=1$
because both $\G1-1$ and $2\G1+1$ are odd.

It is important to emphasize that the   answers obtained in all
the examples considered above (including the latter related to the
divergent series $S_5(k)$ with alternate signs) do  not depend on
the way the items in the entire sum are re-arranged. In fact,
since we always know the exact infinite number of items in the sum
and the order of alternating the signs is clearly defined, we know
also the exact number of positive and negative items in the sum.

Suppose, for instance, that we want to re-arrange the items in the
sum $S_5(\G1)$ in the following way
 \[
 1+1-1+1+1-1+1+1-1+\ldots
\]
However, we know that the sum $S_5(\G1)$ has $\G1$ items and that
grossone is even. This means that in the sum there are
$\frac{\G1}{2}$ positive and $\frac{\G1}{2}$ negative items. As a
result, the re-arrangement considered above can continue only
until the positive items will not finish and then it will be
necessary to continue to add only negative numbers. More
precisely, we have \beq
S_5(\G1)=\underbrace{1+1-1+1+1-1+\ldots+1+1-1}_{\frac{\mbox{\tiny{\G1}}}{2}\mbox{
positive and }\frac{\mbox{\tiny{\G1}}}{4}\mbox{ negative
items}}\hspace{1mm}\underbrace{-1-1-\ldots-1-1-1}_{\frac{\mbox{\tiny{\G1}}}{4}\mbox{
negative items}}=0, \label{Riemann28}
 \eeq
 where the result of the
first part in this re-arrangement is calculated as
$(1+1-1)\cdot\frac{\mbox{\tiny{\G1}}}{4}=\frac{\mbox{\tiny{\G1}}}{4}$
and the result of the second part is equal to
$-\frac{\mbox{\tiny{\G1}}}{4}$.

Note that the record (\ref{Riemann28}) is a re-arrangement of the
sum introduced in (\ref{Riemann27}) where the order of the
alternating of positive and negative items has been defined.
Thanks to this order and to the knowledge of the whole number of
items in the sum (\ref{Riemann27}) we were able to calculate the
number of negative and positive items. Obviously, if we consider
another sum where the order of the alternating of the signs and
the number of the items are given in a different way then the
result can be also different. For instance, if we have the sum
\[
S_6(\G1)=\underbrace{1+1-1+1+1-1+\ldots }_{ \G1 \mbox{ items} }
\]
then we have $\frac{2\mbox{\tiny{\G1}}}{3}$ positive items and
$\frac{\mbox{\tiny{\G1}}}{3}$ negative items that gives the result
$S_6(\G1)=\frac{\mbox{\tiny{\G1}}}{3}$.

Let us consider now the following divergent series
 \[
S_7=  1 -2   + 3 -4 + \ldots
  \]
The corresponding sum $S_7(k)$ can   be easily calculated   as the
difference of two arithmetic progressions after we have fixed the
infinite number of items, $k$, in it. Suppose that $k=\G1$. Then
it follows
 \[
S_7(\G1) = \underbrace{1 -2   + 3 -4 + \ldots - (\G1-2) + (\G1-1)
- \G1}_{ \G1
 \mbox{ items}} =
  \]
  \[
\underbrace{(1     + 3   + 5 + \ldots + (\G1-3) + (\G1-1)
)}_{\frac{\mbox{\tiny{\G1}}}{2}\mbox{ items}} - \underbrace{( 2 +
4 + 6 + \ldots + (\G1-2) + \G1
)}_{\frac{\mbox{\tiny{\G1}}}{2}\mbox{ items}} =
  \]
    \beq
   \frac{(1+(\G1-1)) \G1  }{4} -
\frac{(2+\G1) \G1  }{4} = \frac{ \G1^2 - 2 \G1  - \G1^2}{4}
 = - \frac{\G1}{2}.
 \label{Riemann15}
\eeq Obviously, if we change the number of items, $k$, then, as it
happens in the finite case, the results of summation will also
change. For instance, it follows
 \[  S_7(\G1-1)=\frac{\G1}{2}, \hspace{1cm}
 S_7(\G1+1)=\frac{\G1}{2}+1, \hspace{1cm} S_7(\G1^2)=-\frac{\G1^2}{2}.
\]
In particular, for  $k=\G1-1$ we have
 \[
S_7(\G1-1) = \underbrace{1 -2   + 3 -4 + \ldots - (\G1-2) +
(\G1-1)}_{ \G1 -1\mbox{ items}}   =
  \]
  \[
\underbrace{(1     + 3   + 5 + \ldots + (\G1-3) + (\G1-1)
)}_{\frac{\mbox{\tiny{\G1}}}{2}\mbox{ items}} - \underbrace{( 2
+  4 + 6 + \ldots + (\G1-2)
)}_{\frac{\mbox{\tiny{\G1}}}{2}-1\mbox{ items}} =
  \]
    \[
   \frac{(1+(\G1-1)) \G1  }{4} -
\frac{(2+(\G1-2)) (\G1/2-1)  }{2} = \frac{ \G1^2  }{4} - \left(
\frac{ \G1^2  }{4}  - \frac{\G1}{2}\right)
 =   \frac{\G1}{2}.
 \]
Obviously, we have (cf.  (\ref{Riemann15})) that
 \[
S_7(\G1) =S_7(\G1-1)   -  \G1   = \frac{\G1}{2} -  \G1   = -
\frac{\G1}{2}.
  \]

 Analogously to the passage from series to
sums, the introduction of infinite and infinitesimal numerals
allows us (in fact, it \textit{imposes}) to substitute improper
integrals of various kinds by integrals defined in a more precise
way. For example, let us consider the following improper integral
\[
\int_{0}^{\infty}x^2 dx.
 \]
Now  it is necessary to define its upper infinite limit  of
integration explicitly. Then,   different infinite numbers used
instead of $\infty$ will lead to different results, as it happens
in the finite case. For instance, numbers \G1 and $\G1^2$ give us
two different integrals both assuming infinite but different
values:
 \[
  \int_{0}^{\G1}x^2 dx  = \frac{1}{3}\G1^3,
  \hspace{1cm} \int_{0}^{\G1^2}x^2 dx  = \frac{1}{3}\G1^6.
\]
Moreover, it becomes possible to calculate integrals where both
endpoints of the interval of integration are infinite as in the
following example
 \[
\int_{\G1}^{\G1^2}x^2 dx
  = \frac{1}{3}\G1^6 - \frac{1}{3}\G1^3.
\]

We conclude this section by emphasizing that, as it was with sums,
 it becomes also possible to calculate
  integrals of functions assuming infinite and infinitesimal
  values and the obtained results of integration can be finite, infinite, and
infinitesimal. For instance, in the integral
 \[
\int_{\G1}^{\G1+\G1^{-2}}x^2 dx
  = \frac{1}{3}(\G1^{1}+\G1^{-2})^3 - \frac{1}{3}\G1^3
  =
  1\G1^{0}+1\G1^{-3}+\frac{1}{3}\G1^{-6}
\]
  the result has a finite part and two infinitesimal parts and  in the integral
 \[
\int_{\G1}^{\G1+\G1^{-2}}{x^2-x\,\,\,} dx
  =
  1\G1^{-3}-\frac{1}{2}\G1^{-4}+\frac{1}{3}\G1^{-6}
\]
  the result has  three infinitesimal parts. The last two
  examples illustrate situations when the integrand is infinite

 \[
\int_{\G1}^{\G1^2} \G1 x^2 dx = \G1 \int_{\G1}^{\G1^2} x^2 dx
  = \frac{1}{3}\G1^7 - \frac{1}{3}\G1^4,
\]
and infinitesimal
 \[
\int_{\G1}^{\G1^2} \G1^{-4} x^2 dx  = \G1^{-4} \int_{\G1}^{\G1^2}
x^2 dx
  = \frac{1}{3}\G1^{2}  - \frac{1}{3}\G1^{-1}.
\]

\section{Methodological consequences and analysis of
some \\well known results  related to the zeta function}
\label{s4_riem}

Let us return now to the Riemann zeta function and consider it
together with some Euler's results from the new methodological
positions. The first remark we can make is that the records of the
type
 \beq
 \begin{array}{rl}
f_1(x)= \sum_{i=1}^{\infty} a_i(x), &
\hspace{10mm}f_2(x)=\sum_{-\infty}^{\infty} a_i(x), \\ &\\ f_3(x)=
\int_{a}^{\infty} g(t,x) dt,& \hspace{10mm} f_4(x)=
\int_{-\infty}^{\infty} g(t,x) dt,
\end{array}
 \label{Riemann23}
 \eeq
and similar are   sufficiently precise to define a function only
if one uses the traditional language. As it has been shown in the
previous section,  the accuracy of the answers one can get with
this language is lower with respect to the new numeral system
using \G1.

Thanks to   this numeral system, we know now that the only symbol
$\infty$ and the absence of numerals allowing us to express
different infinite numbers just did not give us the possibility to
distinguish many different functions in each of the objects
present in (\ref{Riemann23}). We can say that we are now
 in the situation similar
to a physicist who first observed a black dot in a weak lens and
then had put a stronger one and has seen that that the black dot
consists of many different dots. When we have changed the lens in
our mathematical microscope we have seen that within each object
of our study viewed as one entity with the old lens we can
distinguish many different separate objects. This means that we
cannot use the traditional language when the accuracy of the
answers to questions made with respect to (\ref{Riemann23}) is
expected to be as high as possible.

Now, when we are aware both of the existence of different infinite
numbers and of the numeral system of Pirah\~{a}, records like
(\ref{Riemann23}) can be written also as
\[
f_1(x)=\sum_{i=1}^{many}a_i(x),
\hspace{10mm}f_2(x)=\sum_{-many}^{many}a_i(x),
\]
\[
 f_3(x)=
\int_{a}^{many} g(t,x) dt, \hspace{10mm} f_4(x)=
\int_{-many}^{many} g(t,x) dt.
 \]
 We emphasize again that records
like (\ref{Riemann23}) are not wrong, they are just less precise
than records that could be written using various infinite
numerals. In the moment when a mathematician  has decided which
numeral system he/she will use in his/her work, he/she has decided
which mathematical objects he/she will be able to observe and to
describe. Thus, if the accuracy of the records (\ref{Riemann23})
is sufficient for the problem under consideration, than everything
is fine. Is this the case of the problem of zeros of the Riemann
zeta function? Let us see.

The fact that within $\infty$ many different infinite numbers can
be distinguished means that records (\ref{Riemann_zeta}) and
(\ref{Riemann_eta}) do not describe two single functions and there
exist \textit{many different} Riemann zeta functions  (and,
obviously, many different Dirichlet eta functions)  and to fix a
concrete one it is necessary to fix the number of items in the
corresponding sum. In order to define a Riemann zeta function (a
Dirichlet eta function), we should choose an infinite number $n$
expressible in a numeral system using grossone (for instance, we
can take the numeral system briefly described in
Section~\ref{s2_riem}) and then write
 \beq
 \zeta(s,n)=\sum_{u=1}^{n}
 \frac{1}{u^s},  \hspace{10mm}
 \eta(s,n)=\sum_{u=1}^{n} \frac{{(-1)}^{u-1}}{u^s}.
\label{Riemann_etazeta}
 \eeq
 As a result, for different
(infinite or finite) values of $n$ we have different functions.
For example,  for $n=\G1/2$ and $n=\G1$ we get
\[
 \zeta(s,\G1/2)=\sum_{u=1}^{\mbox{\tiny{\G1}}/2} \frac{1}{u^s},\hspace{10mm}
 \zeta(s,\G1)=\sum_{u=1}^{\mbox{\tiny{\G1}}} \frac{1}{u^s}
\]
that are two different functions and, analogously,
\[
\eta(s,\G1)=\sum_{u=1}^{\mbox{\tiny{\G1}}}
\frac{{(-1)}^{u-1}}{u^s}, \hspace{10mm}
 \eta(s,\G1/2)=\sum_{u=1}^{\mbox{\tiny{\G1}}/2} \frac{{(-1)}^{u-1}}{u^s},
 \]
are also two different functions.

Note that there obviously exist questions for which both
mathematical languages used in (\ref{Riemann_zeta}) and
(\ref{Riemann_etazeta}) are sufficiently accurate. This can happen
when the property we ask about   holds for all (finite and
infinite) values of $n$ in (\ref{Riemann_etazeta}). For instance,
it is possible easily to answer in the affirmative to the
question: `Is the inequality $ \zeta(1)> -100$ correct?' because
this result holds for all values of $n$
in~(\ref{Riemann_etazeta}).

The above observations mean that   it becomes necessary to
reconsider  classical results concerning properties of the Riemann
zeta function (\ref{Riemann_zeta})  in order to re-write them
(where it is possible) for the form (\ref{Riemann_etazeta}) with
infinite values of $n$. First, it is easy to understand that  the
relation (\ref{Riemann_relation}) can be re-written as
 \beq
 \eta(s,n) = \zeta(s,2k) -  2^{1-s} \zeta(s,k)
    \label{Riemann_relation_new}
 \eeq
for even values of $n=2k$ and as
 \beq
 \eta(s,n) = \zeta(s,2k+1) -  2^{1-s} \zeta(s,k)
    \label{Riemann_relation_new_odd}
 \eeq
 for $n=2k+1$.

Let us look now at the formula (\ref{Riemann_primes}).
 It has been first introduced and studied
by Euler (see \cite{Euler_1,Euler_2,Euler_3,Euler_4}) who proved
for integer values of $s$ the famous identity
 \beq
   \sum_{u=1}^{\infty} \frac{1}{u^s} =
    \prod_{p\,\, \mbox{\tiny{primes}}}  \frac{1}{1-p^{-s}}
     \label{Riemann24}
 \eeq
being one of the strongest sources of the interest to the Riemann
zeta function. This identity has been proved using the traditional
language and, as a consequence, its accuracy depends on the
accuracy of the used language, i.e., the language that cannot
distinguish within $\infty$ different infinite numbers.

If one wants to use the new language and to have a higher accuracy
then it is not possible to write with nonchalance identities
involving infinite sums and/or products of the type
$$\sum_{u=1}^{\infty} a_i =
    \prod_{i=1}^{\infty} b_i.   $$
  It is necessary to take care on the number
of infinite items both in the sum and in the product in the same
way as we do it for the finite number of items. These identities
should be substituted, where this is possible, by the identities
$$\sum_{u=1}^{n} a_i =
    \prod_{i=1}^{k} b_i,$$
    where both $n$ and $k$ are infinite (and probably different)
numbers. Proofs of such identities  should be reconsidered with
the accuracy imposed by the new
    language.

It can easily happen that a result being correct    with the
accuracy of one language is not sufficiently precise when a
language with a higher accuracy is used. We recall that for
Pirah\~{a}
    the records $ 2+1=\mbox{`many'}$
and $2+2= \mbox{`many'}$ are correct. This means that in their
language  $ 2+1=2+2= \mbox{`many'}$. Note that also for us, people
knowing that $2+1=3$ and $2+2=4$, the result of Pirah\~{a} is not
wrong, it is just \textit{inaccurate}. Their mathematical language
cannot express the numbers `three' and `four' but for their
purposes this low accuracy of answers is sufficient.

Let us consider    the Euler product formula (\ref{Riemann24}) and
compare  in this occasion the traditional   language working with
$\infty$ with
  the new one using \G1. To prove
(\ref{Riemann24}) Euler expands each of the factors of the right
hand part of (\ref{Riemann24}) as follows
 \beq
\frac{1}{1-\frac{1}{p^s}}  = 1 + \frac{1}{p^s} + \frac{1}{(p^2)^s}
+ \frac{1}{(p^3)^s} +  \ldots
     \label{Riemann25}
 \eeq
Then he observes that their product is therefore a sum of terms of
the form
 \beq
 \frac{1}{(p^{n_1}_1 \cdot p^{n_2}_2 \cdot \ldots \cdot p^{n_r}_r)^s},
     \label{Riemann26}
 \eeq
where $p_1, p_2, \ldots p_r$ are distinct primes and $n_1, n_2,
\ldots n_r$ are natural numbers. Euler then uses the fundamental
theorem of arithmetic (every integer can be written as a product
of primes) and concludes that the sum of all items
(\ref{Riemann26}) is the left hand part of (\ref{Riemann24}).

In the new language this way of doing is not acceptable because,
in order to start, we should indicate the precise infinite number,
$n$, of items in the sum in the left hand part of
(\ref{Riemann24}). This operation gives us different functions
$\zeta(s,n)$ from (\ref{Riemann_etazeta}). Then (\ref{Riemann25})
should be rewritten by indicating the exact infinite number, $k$,
of items in its right hand part and the result of summing up will
be different with respect to (\ref{Riemann24}). Namely, we have
 \beq
1 + \frac{1}{p^s} + \frac{1}{(p^2)^s} + \frac{1}{(p^3)^s} + \ldots
+ \frac{1}{(p^{k-1})^s} =
\frac{1-\frac{1}{(p^{k})^s}}{1-\frac{1}{p^s}}.
 \label{Riemann39}
 \eeq
The accuracy of the language used by Euler did not allow him
either to observe different infinite values of $n$ and $k$ or to
take into account the infinitesimal value~$\frac{1}{(p^{k})^s}$.

It is necessary to emphasize that the analysis made above does not
mean that (\ref{Riemann24}) is not correct. Euler uses the
language involving the symbol $\infty$ and the accuracy of this
language does not allow him (and actually anyone) to perform a
more precise analysis. Analogously, Pirah\~{a} have $
2+1=\mbox{`many'}$ and $2+2= \mbox{`many'}$ and these results are
correct in their language. Now we have numerals 3 and 4 and are
able to obtain more accurate answers and to observe that
$2+1=3\neq 4=2+2$. Both results, $ 2+1=\mbox{`many'} = 2+2$ and
$2+1=3\neq 4=2+2$, are correct but with different accuracies
determined by the languages used for calculations. Both languages
can be used in dependence on the problem one wishes to deal with.

Another metaphor that can help   is the following. Suppose that we
have measured two distances $A$ and $B$ with the accuracy equal to
1~meter and we have found that both of them are equal to 25
meters. Suppose now that we want to measure them with the accuracy
equal to 1 centimeter. Then, very probably, we shall obtain
something like $A=2487$ centimeters and $B=2538$  centimeters,
i.e., $A \neq B$. Both answers, $A = B$ and $A \neq B$, are
correct but with different accuracies and both of them can be used
successfully in different situations. For instance, if one just
wants to go for a walk, then the accuracy of the answer $A = B$
expressed in meters is sufficient. However, if one needs to
connect some devices with a cable, then a higher accuracy is
required and the answer expressed in centimeters should be used.

We are with the Euler product formula in the same situation. It is
correct with the accuracy of the language using $\infty$ because
this language does not allow one to distinguish different infinite
numbers within $\infty$. In the same time, when a language allows
us to distinguish different infinite and infinitesimal numbers, it
follows from (\ref{Riemann39}) that for any infinite (or finite)
$k$ we have
 \beq 1 + \frac{1}{p^s} +
\frac{1}{(p^2)^s} + \frac{1}{(p^3)^s} + \ldots +
\frac{1}{(p^{k-1})^s} \neq \frac{1}{1-p^{-s}}
 \label{Riemann40}
 \eeq
and, as a result, the following theorem holds.
 \begin{theorem}
For prime $p_i$ and both finite and infinite values of $n$ and $k$
it follows that
 \beq
   \sum_{u=1}^{n} \frac{1}{u^s} \neq
    \prod_{i=1}^{k}  \frac{1}{1-p_i^{-s}}.
     \label{Riemann41}
 \eeq
 \end{theorem}

Thus, the choice of the language fixes the accuracy of the
analysis one can perform. In this paper,  we show that the
traditional language using only the symbol $\infty$ is not
sufficiently precise when one wishes to work with the Riemann zeta
function and the Dirichlet eta function.

 Let us comment upon another famous result of Euler related
to the Riemann zeta function -- his solution to the Basel problem
where he has shown that
 \beq
1 + \frac{1}{2^2} + \frac{1}{3^2} + \frac{1}{4^2} + \ldots =
\frac{\pi^2}{6}.
      \label{Riemann34}
 \eeq
 Let us first briefly present Euler's proof and then comment upon
it. Note that from the point of view of modern mathematics this
proof can be criticized from different points of view. Nowadays
there exist many other proofs considered by mathematicians as more
accurate\footnote{This fact is another manifestation of the
continuous mutation of mathematical languages. The views on
accuracy of proofs have changed since Euler's times and the modern
language is considered more accurate.}.  However, since also in
the modern proofs the symbol $\infty$ and the concept of series
are used, the analysis made below can be applied to the other
proofs, as well.

Euler begins with the standard Taylor expansion of $\sin(x)$,
 \beq
\sin (x) = x - \frac{x^3}{3!} +\frac{x^5}{5!}-\frac{x^7}{7!}+
\ldots ,
     \label{Riemann31}
 \eeq
which converges for all $x$. Euler interprets the left hand side
of (\ref{Riemann31}) as an ``infinite polynomial" $P(x)$ that,
therefore, can be written as product of factors based on its
roots. Since the roots of $\sin (x)$ are $\ldots -3\pi, -2\pi,
-\pi, 0, \pi, 2\pi, 3\pi,\,\, \ldots$ then it follows
\[
P(x)=C x(x^2-\pi^2)(x^2-4\pi^2)(x^2-9\pi^2) \ldots
\]
Euler continues by reminding that $\lim_{x \rightarrow 0}
\frac{\sin (x)}{x}=1$ and
\[
\lim_{x \rightarrow 0} \frac{x(x^2-\pi^2)(x^2-4\pi^2)(x^2-9\pi^2)
\ldots}{x} = ( -\pi^2)( -4\pi^2)( -9\pi^2) \ldots
\]
Thus, $C$ should be the reciprocal of the infinite product on the
right and we obtain
 \beq
\sin (x) = x
(1-\frac{x^2}{\pi^2})(1-\frac{x^2}{2^2\pi^2})(1-\frac{x^2}{3^2\pi^2})
\ldots
     \label{Riemann32}
 \eeq
  Then, as follows from (\ref{Riemann31}) and
(\ref{Riemann32}), Euler has written $P(x)$ in two different ways
and he equates these two records
 \beq
\sin (x) = x - \frac{x^3}{3!} +\frac{x^5}{5!}-\frac{x^7}{7!}+
\ldots = x
(1-\frac{x^2}{\pi^2})(1-\frac{x^2}{2^2\pi^2})(1-\frac{x^2}{3^2\pi^2})
\ldots ,
     \label{Riemann33}
 \eeq
Now Euler equates the coefficients of $x^3$ on both sides of
(\ref{Riemann32}) and  gets first
 \[
-\frac{1}{3} =
-\frac{1}{\pi^2}-\frac{1}{2^2\pi^2}-\frac{1}{3^2\pi^2}
-\frac{1}{4^2\pi^2} \ldots ,
 \]
 and then the final beautiful result (\ref{Riemann34}).

Let us now consider these results using the new approach. We take
the function $\sin (x)$ introduced using the standard
trigonometric reasoning. Note that such a definition just
describes its properties and does not tell us how \textit{to
calculate} $\sin (x)$ precisely at all $x$. When one uses the
trigonometric  definition, values of $\sin (x)$ only at certain
$x$ are known precisely, e.g., $\sin (\frac{\pi}{2})=1$, $\sin
(2\pi)=0$, etc. Notice that in  these records we do not use any
approximation of the number $\pi$. Similarly to $\sin (x)$,  it is
described by its properties and a special numeral, $\pi$, is
introduced to indicate it. Then, those values of $\sin (x)$ that
are not linked to geometric ideas are defined through various
approximations of both $\sin (x)$ and $\pi$.

The new language taken together with the trigonometric definition
of $\sin (x)$   allows us to evaluate $\sin (x)$ precisely not
only at certain finite points but also at certain infinite points.
For instance, it follows $\sin (2\G1 \pi)=0$, $\sin (\G1 \pi+
\frac{\pi}{2})=1$, etc.

Then, the explicit usage of infinite and infinitesimal numbers
imposes us to move from  (\ref{Riemann31}) to the approximation
with the polynomial  $P_1(x,2k+1)$ of the order $2k+1$   where
  \beq
\sin (x) \approx P_1(x,2k+1) =   x - \frac{x^3}{3!}
+\frac{x^5}{5!}-\frac{x^7}{7!}+ \ldots +(-1)^k
\frac{x^{2k+1}}{(2k+1)!}, \hspace{5mm} k=0,1,2, \ldots
\label{Riemann35}
 \eeq
and for different finite or infinite $k$ we get different
approximations. Analogously, the idea used in (\ref{Riemann32})
gives us the second kind of approximation where the polynomial
$P_2(x,2n+1)$ of the order $2n+1$ is used
 \beq
    \begin{array}{rl}
  \sin (x) \approx P_2(x,1) = & x, \\
 \sin (x) \approx P_2(x,2n+1) = &
x(1-\frac{x^2}{\pi^2})(1-\frac{x^2}{2^2\pi^2}) \ldots
(1-\frac{x^2}{n^2\pi^2}), \hspace{2mm}  n= 1,2, \ldots
\end{array}
     \label{Riemann36}
 \eeq
where for different finite or infinite values of $n$ we get
different approximations.

Obviously, when $k\neq n$, it follows $P_1(x,2k+1) \neq
P_2(x,2n+1)$. In the case $k=n$, two polynomials will be equal if
 all their coefficients are equal. Following Euler we
first equate the coefficients of $x^3$ and obtain
 \[
-\frac{1}{3} =
-\frac{1}{\pi^2}-\frac{1}{2^2\pi^2}-\frac{1}{3^2\pi^2}
-\frac{1}{4^2\pi^2}\,\,\, \ldots \,\,\,-\frac{1}{k^2\pi^2},
 \]
From where we get
 \beq
\frac{\pi^2}{6} =1 + \frac{1}{2^2} + \frac{1}{3^2} + \frac{1}{4^2}
+ \ldots + \frac{1}{k^2}.
      \label{Riemann37}
 \eeq
However, such a kind of equating should be done for \textit{all
the coefficients} of $P_1(x,2k+1) $ and $P_2(x,2n+1)$. Suppose
that $k$ in the sum (\ref{Riemann35}) is even, by equating the
coefficients of the highest power $2k+1$ of $x$ in two polynomials
we get that it should be
  \[
 \frac{1}{(2k+1)!}=   \frac{1}{(k!)^2\pi^{2k}},
      \]
from where we have
 \beq
 \frac{\pi^{2k}}{(2k+1)!}=   \frac{1}{(k!)^2}.
      \label{Riemann38}
 \eeq
 Equating coefficients for other powers of $x$ will give us more
different expressions to be satisfied and it is easy to see that
they cannot hold \textit{all together}.

\begin{figure}[t]
  \begin{center}
    \epsfig{ figure = 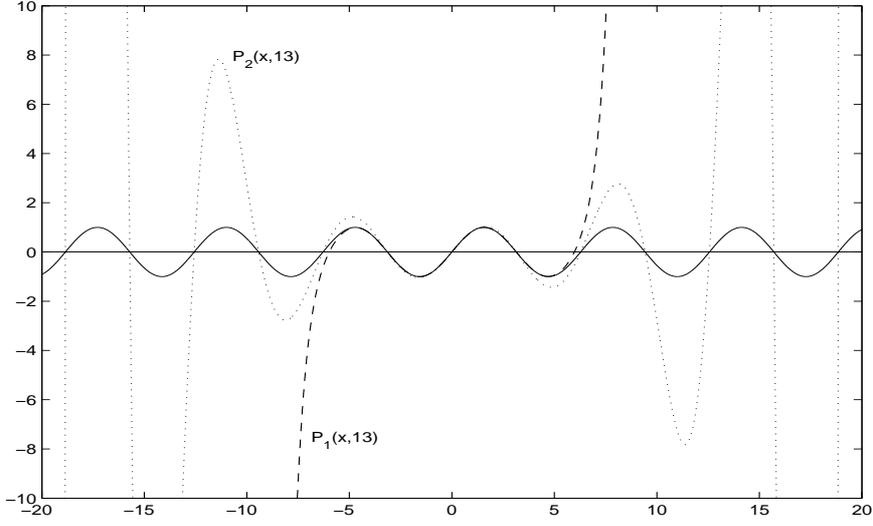, width = 4.5in, height = 2.7in,  silent = yes }
    \caption{Polynomials $P_1(x,13)$ and $P_2(x,13)$ together with
$\sin(x)$. } \label{Riemann_fig}
  \end{center}
\end{figure}

Thus, polynomials $P_1(x,2k+1) $ and $P_2(x,2n+1)$ give us
different approximations of $\sin(x)$. As an example we show in
Fig.~\ref{Riemann_fig} polynomials $P_1(x,13)$ and $P_2(x,13)$
together with $\sin(x)$. Higher is the order of the polynomials
better are the approximations but the differences are always
present for both finite and infinite values of $k$ and $n$.

As it was with the  product formula, Euler  results are correct
with the   accuracy of the traditional mathematical language
working only with the symbol $\infty$. It does not allow one to
distinguish different infinite values of $k$ and $n$ and, as a
consequence, different infinite polynomials $P_1(x,2k+1)$ and
$P_2(x,2n+1)$.

Since all the traditional computations are executed with finite
values of $k$ and $n$ and only the initial finite part of
(\ref{Riemann38}) is used,  it is always possible to choose
sufficiently high finite values of $k$ and $n$ giving so valid
approximations of $\pi$ with a \textit{finite} accuracy. We are
not able by using the traditional mathematical tools to observe
the behavior of   the polynomials and $\sin(x)$ at various points
$x$ at infinity. Thus, the fact of the difference of these objects
in infinity cannot be noticed by finite numeral systems extended
only by the symbol $\infty$. However, it is perfectly visible when
one uses the new numeral system where different   infinite and
infinitesimal numbers can be expressed and distinguished. Thus, we
can conclude this section by the following theorem that has been
just proved.
 \begin{theorem}
For  both finite and infinite values of $n$ and $k$ it follows
that
 \beq
   \sin(x) \neq P_1(x,2k+1) \neq P_2(x,2n+1).
     \label{Riemann42}
 \eeq
 \end{theorem}

\section{Calculating  $\zeta(s,n)$ and $\eta(s,n)$ for infinite
 values of $n$ at $s=0$ and finite points $s= -1, -2, -3,\ldots$}
\label{s6_riem}

In this section, we establish some results for functions
$\zeta(s,n)$ and $\eta(s,n)$ from (\ref{Riemann_etazeta}) using
the approach presented in Section~\ref{s3_riem}. The analysis made
there allows us to calculate $\zeta(s,n)$ and $\eta(s,n)$   for
infinite values of $n$ at $s=0$ and~$s=-1$:
\[
 \zeta(0,n)=  \underbrace{1+1+1+\ldots+1}_{n} = n,
\]
\[
 \zeta(-1,n)=  1+2+3+\ldots+n  = \frac{(n+1)n}{2}.
\]
From where, for instance, at $n=\G1/2$ we obtain
\[
 \zeta(0,\G1/2)=
 \G1/2, \hspace{5mm} \zeta(0,\G1)=
 \G1,
\]
\[
 \zeta(-1,\G1/2)
 =   \frac{(\G1/2+1)\G1}{4}
 \hspace{5mm} \zeta(-1,\G1)
 = \frac{(\G1+1)\G1}{2}.
\]
Analogously, calculation of the   function, $\eta(s,n)$, at the
points $s=0$ and $s=-1$ for infinite values of $n$ also has been
discussed  above. We obtain
\[
\eta(0,n)=\underbrace{1-1+1-1+1-1+1-\ldots}_{n} = \left \{
\begin{array}{ll} 0, &
  \mbox{if  } n=2k,\\
1, &    \mbox{if  } n=2k+1,\\
 \end{array} \right.
\]
for example, $\eta(0,\G1)=0$ and $\eta(0,\G1^2-1)=1$. Analogously,
in (\ref{Riemann15}) we have seen how it is possible to calculate
$\eta(-1,n)$ for infinite values of $n$. For instance, it follows
$\eta(-1,\G1) = - \frac{\G1}{2}$. Note that these results fit well
relation (\ref{Riemann_relation_new}). For example, for $n= \G1$
it follows
\[
 \zeta(0,\G1)-2 \zeta(0,\G1/2)= 0
  = \eta(0,\G1),
\]
\[
 \zeta(-1,\G1)-4\zeta(-1,\G1/2)= -
 \frac{\G1}{2}
  = \eta(-1,\G1).
\]

Let us present now a general method for calculating the values of
$\zeta(s,n)$ and $\eta(s,n)$ for infinite or finite $n$ at $s$
being negative finite integers. For this purpose we use the idea
of Euler to multiply the identity
 \beq
f(x) = \sum_{i=0}^{n} x^i = 1 + x + x^2 + \ldots + x^n,
 \label{Riemann7}
 \eeq
by $x$. By using the differential calculus developed in
\cite{Dif_Calculus} we can apply this idea   without distinction
for both finite and  infinite values of $n$ and $x\neq 1$ that can
be finite, infinite or infinitesimal. By differentiation of
(\ref{Riemann7}) we obtain
 \beq
f'(x) = \sum_{i=0}^{n-1} i x^{i-1} = 1 + 2x + 3x^2 + \ldots + n
x^{n-1},
 \label{Riemann8}
 \eeq
Since   we have   for  finite, infinite or infinitesimal $x\neq 1$
that
 \beq
  f(x) =
\frac{1-x^{n+1}}{1-x}
 \label{Riemann9}
\eeq then by  differentiation of (\ref{Riemann9}) we have that
 \beq
  f'(x) =
\frac{1+nx^{n+1}-(n+1)x^{n}}{(1-x)^2}.
 \label{Riemann10}
\eeq
 Thus, we can conclude that for finite and infinite   $n$ and
for   finite, infinite or infinitesimal $x\neq 1$ it follows
 \beq
\frac{1+nx^{n+1}-(n+1)x^{n}}{(1-x)^2} = \sum_{i=0}^{n-1} i x^{i-1}
= 1 + 2x + 3x^2 + \ldots + n x^{n-1}.
 \label{Riemann11}
\eeq
 By taking $x=-1$ we obtain that the value of $\eta(-1,n)$
 is
 \[
\eta(-1,n) = 1 -2   + 3 -4 + \ldots + n (-1)^{n-1} =
\frac{1+n(-1)^{n+1}-(n+1)(-1)^{n}}{4}.
 \]
 For example, for $n=\G1$ we obtain the result that
 has been already calculated  in (\ref{Riemann15}) in a different way
 \[
\eta(-1,\G1) = 1 -2   + 3 -4 + \ldots + \G1
(-1)^{\mbox{\tiny{\G1}}-1} =
 \]
  \[
\frac{1+\G1(-1)^{\mbox{\tiny{\G1}}+1}-(\G1+1)(-1)^{\mbox{\tiny{\G1}}}}{4}=
\frac{1-\G1 -(\G1+1)}{4} = - \frac{ \G1}{2}.
 \]
If we multiply   successively (\ref{Riemann11}) by $x$ and use
again differentiation of both parts of the obtained equalities
then it becomes possible to obtain the values of  $\eta(s,n)$ for
other finite integer negative points $s$. Thus, in order to
obtain, for instance, formulae for $s= -2, -3, -4$, we proceed as
follows
 \beq
1 + 2^2x + 3^2x^2 + \ldots + n^2 x^{n-1} =
\frac{-n^{2}x^{n+2}+(2n^2+2n-1)x^{n+1}-(n+1)^{2}x^{n}+x+1}{(1-x)^3},
 \label{Riemann12}
\eeq
\[
1 + 2^3x + 3^3x^2 + 4^3x^3 \ldots + n^3 x^{n-1} =
(n^{3}x^{n+3}-(3n^3+3n^2-3n+1)x^{n+2}+
\]
\beq (3n^3+6n^2-4)x^{n+1}-(n+1)^{3}x^{n}+x^2+4x+1)(1-x)^{-4},
 \label{Riemann13}
\eeq
\[
1 + 2^4x + 3^4x^2 + 4^4x^3 \ldots + n^4 x^{n-1} =
(n^{4}x^{n+4}-(4n^4+4n^3-6n^2+4n-1)x^{n+3}+
\]
\[
(6n^4+12n^3-6n^2-12n+11)x^{n+2}- (4n^4+12n^3+6n^2-12n-11)x^{n+1}+
\]
\beq
 (n+1)^{4}x^{n}-x^3-11x^2-11x-1)(x-1)^{-5}.
 \label{Riemann14}
\eeq
 By taking $x=-1$ in (\ref{Riemann12}), (\ref{Riemann13}), and
(\ref{Riemann14}) we obtain that
 \[
\eta(-2,n)=
2^{-3}(-n^{2}(-1)^{n+2}+(2n^2+2n-1)(-1)^{n+1}-(n+1)^{2}(-1)^{n}),
\]
\[
\eta(-3,n)= 2^{-4}(n^{3}(-1)^{n+3}-(3n^3+3n^2-3n+1)(-1)^{n+2}+
\]
\[
(3n^3+6n^2-4)(-1)^{n+1}-(n+1)^{3}(-1)^{n}-2),
 \]
\[
\eta(-4,n)= -2^{-5}
(n^{4}(-1)^{n+4}-(4n^4+4n^3-6n^2+4n-1)(-1)^{n+3}+
\]
\[
(6n^4+12n^3-6n^2-12n+11)(-1)^{n+2}- (4n^4+12n^3+6n^2
\]
\[
 -12n-11)(-1)^{n+1}+(n+1)^{4}(-1)^{n}).
 \]
 Then, by taking different infinite values of $n$ we are able to
 calculate the respective values of the functions. For example,
 for $n=\G1$ it follows
 \beq
\eta(-2,\G1)= 2^{-3}(-\G1^{2}-2\G1^2-2\G1+1-(\G1+1)^{2})=
-0.5\G1(\G1 +1),
 \label{Riemann20} \eeq
\[
\eta(-3,\G1)= 2^{-4}(-\G1^{3} -(3\G1^3+3\G1^2-3\G1+1) -
\]
\beq
 (3\G1^3+6\G1^2-4)
-(\G1+1)^{3}-2)= -0.5\G1^2(\G1 +3),
 \label{Riemann21}
\eeq
\[
\eta(-4,\G1)= -2^{-5} (\G1^{4}+ 4\G1^4+4\G1^3-6\G1^2+4\G1-1
+6\G1^4+12\G1^3-
\]
\[
 6\G1^2-
 12\G1+11 +  4\G1^4+12\G1^3+6\G1^2
 -12\G1-11  +(\G1+1)^{4})=
 \]
\beq -0.5\G1 (\G1 +1)(\G1^2+\G1-1).
 \label{Riemann22}
\eeq

In order to use the same technique for calculating   $\zeta(s,n)$
 for infinite or finite $n$ at negative
finite integers $s<-1$   it would be necessary to evaluate
(\ref{Riemann12}), (\ref{Riemann13}), and (\ref{Riemann14}) at the
point $x=1$. In traditional mathematics this is impossible whereas
the new approach (see \cite{Korea,Dif_Calculus}) allows us to
execute the required evaluations by using the following method.

If we put $x=1$ at the left-hand parts of (\ref{Riemann12}),
(\ref{Riemann13}), and (\ref{Riemann14}) then we see that we have
there infinite sums of positive integers. Thus, these sums should
be equal to some infinite positive integers. Since in the
right-hand parts of these equalities it is not possible to use
$x=1$, we introduce an infinitesimal perturbation,
$\G1^{-\alpha}$, $\alpha>0$, and calculate the the right-hand
parts  at the point $x=1+\G1^{-\alpha}$ that is infinitesimally
close to $x=1$. Then, in the obtained result, we separate the
contribution of the perturbation that can be kept infinitesimal by
the choice of $\alpha$, from the contribution of the point $x=1$
that should be equal, as we have established,  to an infinite
integer.

In order to present the method, let us execute calculations for
$s= -2$, results for other values of $s$ are obtained by a
complete analogy. We indicate the right-hand part of
(\ref{Riemann12}) as $f(x,n)$ and, by using the usual notation,
${n \choose k}$, for binomial coefficients, proceed as follows
 \[
f(1+\G1^{-\alpha},n)= -\G1^{3\alpha}
\Big(-n^{2}\Big[1+(n+2)\G1^{-\alpha}+
\frac{1}{2}(n+1)(n+2)\G1^{-2\alpha}+
 \]
\[
  \frac{1}{6}n(n+1)(n+2)\G1^{-3\alpha}+
  {\textstyle{n+2 \choose 4}}
  \G1^{-4\alpha} + \ldots
+ \G1^{-(n+2)\alpha}\Big]+
 \]
\[
(2n^2+2n-1) \Big[ 1+(n+1)\G1^{-\alpha}+
\frac{1}{2}n(n+1)\G1^{-2\alpha}+
\frac{1}{6}(n-1)n(n+1)\G1^{-3\alpha}+
\]
\[
{\textstyle{n+1 \choose 4}}
  \G1^{-4\alpha} +\ldots + \G1^{-(n+1)\alpha}\Big]-
(n+1)^{2}\Big[1+ n \G1^{-\alpha}+ \frac{1}{2}(n-1)n
\G1^{-2\alpha}+
 \]
\[
\frac{1}{6}(n-2)(n-1)n\G1^{-3\alpha}+ {\textstyle{n \choose 4}}
  \G1^{-4\alpha} + \ldots +
\G1^{- n\alpha}\Big]+\G1^{-\alpha}+2\Big).
\]
By collecting the terms of grossone we then obtain
\[
f(1+\G1^{-\alpha},n)= \frac{1}{6}n(n+1)(2n+1) +
 \]
\beq
 \Big(n^{2} {\textstyle{n+2 \choose 4}} -(2n^2+2n-1){\textstyle{n+1 \choose
4}}+(n+1)^{2}{\textstyle{n \choose 4}}\Big)\G1^{-\alpha}+ \ldots +
n^{2}\G1^{-(n+2)\alpha}.
  \label{Riemann16}
\eeq
 As it can be seen from (\ref{Riemann16}), for any finite or infinite
value of $n$ there always can be chosen a  number $\alpha>0$ such
that the contribution of the added infinitesimal $\G1^{-\alpha}$
in $f(1+\G1^{-\alpha},n)$ is a sum of infinitesimals (see the
second line of   (\ref{Riemann16})). Due to the representation
(\ref{3.12}), (\ref{3.13}), this contribution can be easily
separated from the integer finite or infinite part represented by
the first line of (\ref{Riemann16}).

Thus, we have obtained that for finite and infinite values of $n$
it follows
 \beq
  \zeta(-2,n)=\frac{1}{6}n(n+1)(2n+1).
 \label{Riemann17}
\eeq
Analogously, by applying the same procedure again we can
obtain the formulae for $\zeta(s,n)$ for other finite integer
negative points $s$. For instance,
 \beq
 \zeta(-3,n)=
\frac{1}{4}n^2(n+1)^2,
 \label{Riemann18}
\eeq
 \beq
\zeta(-4,n)= \frac{1}{30}n(n+1)(2n+1)(3n^2+3n-1).
 \label{Riemann19}
\eeq

Note that the obtained results fit perfectly both well-known
formulae for finite values of $n$ (see \cite{Beyer}) and relation
(\ref{Riemann_relation_new}). For example, by using
(\ref{Riemann_relation_new}) and the obtained formulae
(\ref{Riemann17}), (\ref{Riemann18}), and (\ref{Riemann19}) with
$n= \G1$ and $n= \G1/2$ we obtain (cf. (\ref{Riemann20}),
(\ref{Riemann21}), and (\ref{Riemann22})) that
\[
 \zeta(-2,\G1)-8 \zeta(-2,\G1/2)=
 \frac{1}{6}\G1(\G1+1)(2\G1+1)
 -
 \frac{1}{3}\G1(\G1+2)(\G1+1)=
\]
\[
 -0.5\G1(\G1 +1) = \eta(-2,\G1),
\]
\[
 \zeta(-3,\G1)-16\zeta(-3,\G1/2)=
 \frac{1}{4}\G1^2(\G1+1)^2-
  \G1^2(0.5\G1+1)^2
  =
\]
\[
 -0.5\G1^2(\G1 +3) =  \eta(-3,\G1),
\]
\[
 \zeta(-4,\G1)-32\zeta(-4,\G1/2)=
 \frac{1}{30}\G1(\G1+1)(2\G1+1)(3\G1^2+3\G1-1)-
 \]
\[
\frac{1}{15}\G1(\G1+2)(\G1+1)(3\G1^2+6\G1-4)= -0.5\G1 (\G1
+1)(\G1^2+\G1-1)
  = \eta(-4,\G1).
\]

\section{Conclusion}
\label{s7_riem}

There has been shown in this paper that, as it happens in Physics,
in Mathematics the instruments used to observe mathematical
objects have their own accuracy and the chosen accuracy bounds and
influences results of the observation.  Numeral systems used to
represent numbers are among the instruments of mathematicians.
When a mathematician chooses a mathematical language (an
instrument), in this moment he/she chooses both a set of numbers
that can be observed through the numerals available in the chosen
numeral system and the accuracy of results that can be obtained.
In the cases where two languages having different accuracies can
be applied, it does not usually make sense to mix the languages,
i.e., to compose mathematical expressions using symbols from both
languages, because the result of such a mixing either has no any
sense or has the lower of the two accuracies.

Two mathematical languages have been studied in the paper: (i)~the
traditional mathematical language using the symbol $\infty$;
(ii)~the new language allowing one to represent different infinite
and infinitesimal numbers. These languages have been applied and
compared in several contexts related to the Riemann zeta function
and the Dirichlet eta function. It has been discovered that the
situations that can be illustrated by  the following metaphor can
take place.

Suppose that we have measured two distances $A$ and $B$ with the
accuracy equal to 1~meter and we have found that both of them are
equal to 25 meters. Suppose now that we want to measure them with
the accuracy equal to 1 centimeter. Then, very probably, we shall
obtain something like $A=2487$ centimeters and $B=2538$
centimeters, i.e., $A \neq B$. Both answers, $A = B$ and $A \neq
B$, are correct but with different accuracies and both of them can
be used successfully in different situations. For instance, if one
just wants to go for a walk, then the accuracy of the answer $A =
B$ expressed in meters is sufficient. However, if one needs to
connect some devices with a cable, then a higher accuracy is
required and the answer expressed in centimeters should be used.

The analysis done in the paper shows that the traditional
mathematical language using the symbol $\infty$ very often does
not possess a sufficiently high accuracy when one deals with
problems having their interesting properties at infinity. For
instance, it does not allow us to distinguish within the record
$f(x)=\sum_{i=1}^{\infty}a_i(x)$ different functions
$f_n(x)=\sum_{i=1}^{n}a_i(x)$ emerging for different infinite $n$.
However,  functions $f_n(x)$ become visible if one uses the more
powerful numeral systems allowing us to write down different
infinite (and infinitesimal) numbers. Then,  it follows that if
$a_{n+1}(x)\neq 0$ functions $f_n(x)$ and $f_{n+1}(x)$  are
different  and $f_n(x) \neq f_{n+1}(x)$. When $a_{n+1}(x)$ is an
infinitesimal number then the difference $f_n(x) - f_{n+1}(x)$ is
also infinitesimal, i.e., invisible if one uses the  traditional
mathematical language  but perfectly observable through the new
numeral system.

In particular, this means that it is impossible to answer to the
Riemann Hypothesis because it asks a question about the behavior
of something that there was supposed to be a function but, as we
can see using the new numeral system distinguishing different
infinite numbers,   is not a function but many different
functions. Moreover, the analysis of the Riemann zeta function
made traditionally does not consider various infinite and
infinitesimal numbers that are crucial not only in the context of
the Hypothesis but even with respect to the definition of the
function itself.

Of course, there remains a possibility that the question asked in
the Hypothesis has the same answer for any infinite $n$.
Unfortunately,  this is impossible due to the analysis made above
and    the analysis of the partial zeta functions made in
\cite{Balazard,Borwein,Spira}.

We conclude this paper by the following remark. It is well known
that the Riemann zeta function appears in many different
mathematical contexts and there exist many different equivalent
formulations of the Riemann hypothesis. In  view of the above
analysis, this can be explained again by the accuracy of the
traditional mathematical language. Such manifestations indicate
situations where its accuracy is not sufficient. By returning to
the metaphor with the microscope, we can say that out traditional
lens is too weak to distinguish different mathematical objects
observed in these situations. Thus,   these mathematical contexts
can be viewed as very promising for obtaining new results by
applying the new numeral system that does not use the symbol
$\infty$ and allows one to distinguish and to treat
 numerically various infinite and infinitesimal numbers.

\bibliographystyle{plain}
\bibliography{XBib_accuracy}

\end{document}